\documentclass[review,onefignum,onetabnum]{siamart250211}

\usepackage{lipsum}
\usepackage{amsfonts}
\usepackage{graphicx}
\usepackage{epstopdf}
\usepackage{algorithmic}
\usepackage{amssymb}
\usepackage[english]{babel}
\usepackage{array}
\usepackage{adjustbox}
\usepackage{relsize}
\usepackage{spverbatim}
\usepackage{listings}
\usepackage{mathrsfs}
\usepackage{cleveref}
\usepackage{derivative}
\usepackage{mathtools}
\usepackage{float}
\usepackage[section]{placeins}
\usepackage[numbers,sort]{natbib}

\DeclareMathAlphabet{\mathpzc}{OT1}{pzc}{m}{it}

\def\Frac#1#2{\displaystyle{\frac{#1}{#2}}}

\newcolumntype{L}{>{$}l<{$}}

\ifpdf
  \DeclareGraphicsExtensions{.eps,.pdf,.png,.jpg}
\else
  \DeclareGraphicsExtensions{.eps}
\fi

\newsiamremark{hypothesis}{Hypothesis}
\crefname{hypothesis}{Hypothesis}{Hypotheses}
\newsiamthm{claim}{Claim}

\nolinenumbers

\headers{Zeros of Reverse Generalized Bessel Polynomials}{T. M. Dunster, A. Gil, D. Ruiz-Antolin and J. Segura}

\title{Uniform Asymptotic approximation and numerical evaluation of the Reverse Generalized Bessel Polynomial zeros}

\author{T. M. Dunster\thanks{Department of Mathematics and Statistics, San Diego State University, 5500 Campanile Drive, San Diego, CA 92182, USA. 
  (\email{mdunster@sdsu.edu}, \url{https://tmdunster.sdsu.edu}).}
\and A. Gil\thanks{Departamento de Matem\'atica Aplicada y CC. de la Computaci\'on, ETSI Caminos, Universidad de Cantabria, 39005-Santander, Spain. 
  (\email{amparo.gil@unican.es}).}
  \and D. Ruiz-Antolin\thanks{Departamento de Matem\'atica Aplicada y CC. de la Computaci\'on, ETSI Caminos, Universidad de Cantabria, 39005-Santander, Spain. 
  (\email{diego.ruizantolin@unican.es}).}
\and J. Segura\thanks{Departamento de Matem\'aticas, Estad\'{\i}stica y Computaci\'on, Facultad de Ciencias, Universidad de Cantabria, 39005-Santander, Spain. 
  (\email{javier.segura@unican.es}).}}
\usepackage{amsopn}

\makeatletter
\newcommand*{\addFileDependency}[1]{
  \typeout{(#1)}
  \@addtofilelist{#1}
  \IfFileExists{#1}{}{\typeout{No file #1.}}
}
\makeatother

\nolinenumbers
\ifpdf
\hypersetup{
 pdftitle={Uniform Asymptotic approximation and numerical evaluation of the Reverse Generalized Bessel Polynomial zeros},
 pdfauthor={T. M. Dunster, A. Gil, D. Ruiz-Antolin and J. Segura}
}
\fi
\begin{document}
\maketitle
\begin{abstract}
Uniform asymptotic expansions are derived for the zeros of the reverse generalized Bessel polynomials of large degree $n$ and real parameter $a$. It is assumed that $-\Delta_{1} n+\frac{3}{2} \leq a \leq \Delta_{2} n$ for fixed arbitrary $\Delta_{1} \in (0,1)$ and bounded positive $\Delta_{2}$. For this parameter range at most one of the zeros is real, with the rest being complex conjugates. The new expansions are uniformly valid for all the zeros, and are shown to be highly accurate for moderate or large values of $n$. They are consequently used as initial values in a very efficient numerical algorithm designed to obtain the remaining complex zeros using Taylor series. 
\end{abstract}


\section{Introduction}
\label{sec:Introduction}

The generalized Bessel polynomials are defined by
\begin{equation} 
\label{eq01}
y_{n}(z;a)=\sum_{k=0}^{n}\binom{n}{k}(n+a-1)_{k}\left(\tfrac12 z\right)^{k},
\end{equation}
where $(\alpha)_{k}=\Gamma(\alpha+k)/\Gamma(\alpha)$ is Pochhammer’s symbol. Their zeros, which are generally complex-valued, arise in a number of applications in applied mathematics \cite{Kong:2024:EIG} and engineering: see, for example, \cite{Johnson:1976:FIL,Fila:2015:AFI,Carimalo:2018:MFG,Martinez:1977:TFG}. Previous works addressing the approximation or computation of these zeros include \cite{Carpenter:1992:GBP,Pasquini:2000:GBP}. In addition, \cite{Segura:2013:CCZ} presents a general method for computing complex zeros of special functions, one of the cases 
considered being the zeros of reverse generalized Bessel polynomials. For an investigation into the domains in which the zeros lie, see \cite{deBruin:1981:ZGBI,deBruin:1981:ZGBII}.

In this paper, we derive uniform asymptotic expansions for these zeros as $n \to \infty$, which are uniformly valid for \emph{all} the zeros, and for $-\Delta_{1} n+\frac{3}{2} \leq a \leq \Delta_{2} n$ for fixed arbitrary $\Delta_{1} \in (0,1)$ and bounded positive $\Delta_{2}$. These expansions are considerably more powerful than existing results, since they are uniformly valid for $n$ large and $|a|$ small or large. These are derived in \cref{sec:zeros}, employing certain coefficients presented in \cref{sec:LG}, and are shown to be highly accurate for moderate or large values of $n$. Moreover, they can be used as starting values in a new highly efficient numerical algorithm which we develop in \cref{sec:algorithm}, which is designed to obtain the remaining complex zeros using Taylor series. 

It is important to mention that the matrix method presented in 
\cite{Pasquini:2000:GBP} for calculating the zeros of generalized Bessel polynomials requires good initial approximations and is computationally expensive. In contrast, our algorithm avoids these drawbacks, even when many zeros are required. Moreover, using Taylor series avoids explicit evaluation of the function: since we only need the zeros, the overall normalization is irrelevant, simplifying the computation.

In our analysis, following \cite{Dunster:2001:GBP,Dunster:2025:SAR,Dunster:2021:CPB}, we find it more convenient to consider the reverse Bessel polynomials, which are given by
\begin{equation} 
\label{eq02}
\theta_{n}(z;a)=z^{n}y_{n}(z^{-1};a).
\end{equation}
Then, from \cite[Eqs.~(2.3) and (2.4)]{Dunster:2001:GBP}, define the scaled function
\begin{equation} 
\label{eq03}
w_{n}^{(0)}(z;a)=2^{-n-a+1}z^{1-n-a/2} 
e^{-z} \theta_{n}(z;a),
\end{equation}
which satisfies the differential equation
\begin{equation} 
\label{eq03a}
\frac{d^{2}w}{dz^{2}}=\left\{1+\frac{a-2}{z}
+\frac{(2n+a)(2n+a-2)}{4z^{2}}\right\}w.
\end{equation}
This has a regular singularity at $z = 0$ and an irregular singularity at infinity. The significance of $w_{n}^{(0)}(z;a)$ is that it is the solution that is recessive at infinity in the right half-plane, since
\begin{equation} 
\label{eq04}
w_{n}^{(0)}(z;a)=2^{-n-a+1}z^{1-a/2} e^{-z}
\left\{1+\mathcal{O}(z^{-1})\right\}
\quad (z \to \infty).
\end{equation}

As in \cite{Dunster:2025:SAR}, we also use two numerically satisfactory companion solutions of \cref{eq03a}, namely
\begin{equation} 
\label{eq05}
w_{n}^{(1)}(z;a)=(-1)^{n+1}z^{n+a/2}e^{-z}
V(n+a-1,2n+a,2z),
\end{equation}
and
\begin{equation} 
\label{eq06}
w_{n}^{(-1)}(z;a)=z^{n+a/2}e^{-z}
\mathbf{M}(n+a-1,2n+a,2z),
\end{equation}
where $V(a,b,z)$ and $\mathbf{M}(a,b,z)$ are certain confluent hypergeometric functions (see \cite[pp.~255--256]{Olver:1997:ASF}). $w_{n}^{(1)}(z;a)$ and $w_{n}^{(-1)}(z;a)$ are recessive at $z=\infty$ in the left half-plane $|\arg(-z)|<\tfrac12\pi$, and at $z=0$, respectively. This follows from the limiting behavior
\begin{equation} 
\label{eq07}
w_{n}^{(1)}(z;a)=2^{-n-1}z^{(a/2)-1} e^{z}
\left\{1+\mathcal{O}(z^{-1})\right\}
\quad \left(z \to \infty,\, |\arg(-z)|
\leq \tfrac32 \pi - \delta\right),
\end{equation}
where $\delta$ is an arbitrary small positive constant, and 
\begin{equation} 
\label{eq08}
w_{n}^{(-1)}(z;a)=\frac{z^{n+a/2}}{\Gamma(2n+a)}
\left\{1+\mathcal{O}(z)\right\}
\quad (z \to 0).
\end{equation}

\section{Liouville-Green coefficients}
\label{sec:LG}
We define
\begin{equation} 
\label{eq09}
u= n+\frac12, \, \alpha = \frac{a - 2}{u}.
\end{equation}
Then the differential equation \cref{eq03a} can be rewritten in the form
\begin{equation} 
\label{eq10}
\frac{d^2 w}{dz^2} 
= \left\{ u^2 f(\alpha, z) + g(z) \right\} w,
\end{equation}
where
\begin{equation} 
\label{eq11}
f(\alpha, z) =
\frac{\left(z + \tfrac{1}{2} \alpha\right)^2 +1+\alpha}{z^{2}},
\quad
g(z) = -\frac{1}{4z^2}.
\end{equation}
On factoring we note that
\begin{equation} 
\label{eq12}
f(\alpha, z) = \frac{(z - z_1)(z - z_2)}{z^2},
\end{equation}
where
\begin{equation} 
\label{eq13}
z_{1,2}(\alpha) 
= \pm i \sigma - \tfrac{1}{2} \alpha,
\end{equation}
in which
\begin{equation} 
\label{eq14}
\sigma=\sqrt{1+\alpha}.
\end{equation}
Thus, for large $u$, \cref{eq10} has turning points at $z=z_{1,2}$. Following \cite{Dunster:2025:SAR}, we assume
\begin{equation} 
\label{eq15}
-1<-1+\delta \leq \alpha \leq \alpha_{1} < \infty,
\end{equation}
and, as such, the two turning points are bounded complex conjugates, bounded away from each other and from the pole at $z=0$.

In our expansions for the zeros, we use a Liouville-Green (LG) variable $\xi$, along with the Liouville variable $\zeta$, which appears in turning point expansions (see \cite[Chaps. 10 and 11]{Olver:1997:ASF}). These are given in the present case by
\begin{multline}
\label{eq16}
\frac{2}{3}\zeta^{3/2}=\xi = \int_{z_{1}(\alpha)}^{z} f^{1/2}(\alpha, t)\, dt
\\
=Z-\left(1+\tfrac12 \alpha\right)
\ln\left\{\frac{4Z+2\alpha(Z+z+2)+4+\alpha^{2}}{z} \right\} 
\\
+\tfrac12 \alpha\ln(2Z+2z+\alpha)
+\tfrac{1}{2}\ln(1+\alpha)
+\left(2+\tfrac12 \alpha\right)\ln(2)
-\tfrac{1}{2}(1+\alpha)\pi i,
\end{multline}
where
\begin{equation} 
\label{eq17}
Z = \left\{(z - z_1)(z - z_2)\right\}^{1/2}
= \left\{\left(z + \tfrac{1}{2} \alpha\right)^2 
+1+\alpha \right\}^{1/2}.
\end{equation}

The branch of the square root in \cref{eq17} is chosen so that $Z > 0$ for $z > 0$, $Z < 0$ for $z < 0$, with $Z$ being continuous throughout the upper half of the complex $z$ plane, except along a branch cut connecting $z = 0$ to the turning point $z_1$, where the imaginary part of $\xi$ vanishes. This cut traces the so-called anti-Stokes line. With this choice, we have $Z \sim z$ as $z \to \infty$ in the upper half-plane ($\Im(z) \geq 0$). Furthermore, principal branches are used for the logarithmic terms appearing in \cref{eq16}. A detailed description of the associated conformal mapping is provided in \cite{Dunster:2025:SAR}.

The following coefficients appearing in LG expansions were constructed in \cite{Dunster:2025:SAR}, and we shall use them here. First, let $\phi \in \mathbb{C}$ be defined by
\begin{equation} 
\label{eq18}
\sin(\phi)
=\frac{\sigma}{Z},
\end{equation}
and hence from \cref{eq17},
\begin{equation} 
\label{eq19}
\cos(\phi)
=\frac{z+\tfrac{1}{2} \alpha}{Z}.
\end{equation}
The functions $\sin(\phi)$ and $\cos(\phi)$ are both positive when $z > 0$, and they extend continuously throughout the upper half-plane $\Re(z) \geq 0$, except along the branch cut extending from $z = 0$ to the turning point $z_1$ along the anti-Stokes line as described above. In particular, $\cos(\phi)$ is positive for all $z$ in the interval $(-\infty, -\tfrac{1}{2}\alpha)$, and approaches 1 as $z \to \infty$ along any ray in the upper half-plane. Also observe that, by combining \cref{eq18} and \cref{eq19}, one obtains the relation
\begin{equation} 
\label{eq20}
z=\sigma\,\cot(\phi)
-\tfrac{1}{2}\alpha.
\end{equation}

With these definitions, the LG coefficients that we shall use are given by
\begin{equation} 
\label{eq21}
\mathrm{E}_{1}(\alpha,\phi)
=\frac{\sin(\phi)\left\{5\cos^{2}(\phi)-2\right\}}
{24\,\sigma}+\frac{\alpha\left\{\cos(\phi)
\left(5\cos^{2}(\phi)-6\right)+1 \right\}}{48(1+\alpha)},
\end{equation}
\begin{multline} 
\label{eq22}
\mathrm{E}_{2}(\alpha,\phi)
=\frac{\alpha\cos(\phi)\sin^{3}(\phi)\left\{3-5\cos^{2}(\phi)\right\}}
{16\,(1+\alpha)^{3/2}}
\\
+\frac{\sin^{2}(\phi)}{64(1+\alpha)^{2}}
\left\{ 5\left(4-\alpha^{2}+4\alpha\right)\cos^{4}(\phi)
+(7\alpha^2 - 16\alpha - 16)\cos^{2}(\phi)
-2\alpha^{2} \right\},
\end{multline}
and for $s=2,3,4,\ldots$
\begin{equation} 
\label{eq23}
\mathrm{E}_{s+1}(\alpha,\phi) =
G(\alpha,\phi)
\frac{\partial \mathrm{E}_{s}(\alpha,\phi)}
{\partial \phi}
+\int_{0}^{\phi} 
G(\alpha,\varphi)\sum\limits_{j=1}^{s-1}
\frac{\partial \mathrm{E}_{j}(\alpha,\varphi)}{\partial \varphi}
\frac{\partial \mathrm{E}_{s-j}(\alpha,\varphi)}{\partial \varphi} d\varphi,
\end{equation}
where
\begin{equation} 
\label{eq24}
G(\alpha,\phi)
=-\frac12 \frac{d\phi}{d\xi}
=\frac{\cos(\phi)\sin^{2}(\phi)}
{2 \sigma}
-\frac{\alpha \sin^{3}(\phi)}{4(1+\alpha)}.
\end{equation}
The lower integration limits in \cref{eq23} are chosen for convenience so that $\mathrm{E}_{s}(\alpha,0)=0$, which means they vanish as $z \to \infty$.

Next, let
\begin{equation}
\label{eq25}
a_{1}=a_{2}=\tfrac{5}{72}, \;
\tilde{a}_{1}=\tilde{a}_{2}=-\tfrac{7}{72},
\end{equation}
with subsequent terms $a_{s}$ and $\tilde{{a}}_{s}$ satisfying the same recursion formula, viz.
\begin{equation}
\label{eq26}
a_{s+1}=\frac{1}{2}\left(s+1\right) a_{s}+\frac{1}{2}
\sum\limits_{j=1}^{s-1}{a_{j}a_{s-j}} 
\quad (s=2,3,4,\ldots).
\end{equation}
Further, define
\begin{equation} 
\label{eq27}
\tilde{\mathcal{E}}_{s}(a,z)
=\mathrm{E}_{s}(\alpha,\phi)
+(-1)^{s}\frac{\tilde{a}_{s}}{s\xi^{s}},
\end{equation}
and
\begin{equation} 
\label{eq28}
\mathcal{E}_{s}(a,z)
=\mathrm{E}_{s}(\alpha,\phi)
+(-1)^{s}\frac{a_{s}}{s\xi^{s}}.
\end{equation}
Also, let the sequence $d_{2s+1}(\alpha)$ ($s=0,1,2,\ldots$) be given by
\begin{multline} 
\label{eq29}
\frac12 \biggl[u\alpha(\ln(u)-1)+u(1+\alpha)\ln(1+\alpha) 
\biggr.
\\ 
\left.
+\ln\left\{\Gamma\left(u+\frac12\right)\right\}
-\ln\left\{\Gamma\left(u+u\alpha+\frac12\right)\right\}\right]
\sim \sum_{s=0}^{\infty}\frac{d_{2s+1}(\alpha)}{u^{2s+1}}
\quad (u \to \infty).
\end{multline}
The first four terms are given by
\begin{equation} 
\label{eq30}
d_{1}(\alpha)=-\frac{\alpha}{48(1+\alpha)},
\end{equation}
\begin{equation} 
\label{eq31}
d_{3}(\alpha)=\frac{7\alpha\left(3+3\alpha
+\alpha^{2}\right)}{5760(1+\alpha)^{3}},
\end{equation}
\begin{equation} 
\label{eq32}
d_{5}(\alpha)=-\frac{31\alpha\left(5+10\alpha
+10\alpha^{2}+5\alpha^{3}+\alpha^{4}\right)}
{80640(1+\alpha)^{5}},
\end{equation}
and
\begin{equation} 
\label{eq33}
d_{7}(\alpha)
=\frac{127\alpha\left(7+21\alpha+35\alpha^{2}
+35\alpha^{3}+21\alpha^{4}+7\alpha^{5}+\alpha^{6}\right)}
{430080(1+\alpha)^{7}}.
\end{equation}

Then from \cite[Thm. 3.2]{Dunster:2025:SAR}
\begin{multline} 
\label{eq34}
\theta_{n}(uz;a)
=\frac{u^{1/6}}{e^{(u+a+2)\pi i/2}}
\left\{\frac{2^{a} \pi n! }{\Gamma(n+a-1)}
\right\}^{1/2}\left\{\frac{\zeta}{f(a,z)}\right\}^{1/4}
(uz)^{n+\frac{1}{2}a-1} e^{uz} 
\\ \times
\left\{\mathrm{Ai}\left(u^{2/3}\zeta\right)A(u,a,z)
+\mathrm{Ai}'\left(u^{2/3}\zeta \right)B(u,a,z)\right\},
\end{multline}
where
\begin{equation} 
\label{eq35}
A(u,a,z) \sim 
\exp \left\{ \sum\limits_{s=1}^{\infty}
\frac{\tilde{\mathcal{E}}_{2s}(a,z) }{u^{2s}}\right\} 
\cosh \left\{ \sum\limits_{s=0}^{\infty}
\frac{\tilde{\mathcal{E}}_{2s+1}(a,z)
+d_{2s+1}(\alpha)}
{u^{2s+1}}\right\},
\end{equation}
\begin{equation} 
\label{eq36}
B(u,a,z) \sim \frac{1}{u^{1/3}\zeta^{1/2}}
\exp \left\{ \sum\limits_{s=1}^{\infty}
\frac{\mathcal{E}_{2s}(a,z) }{u^{2s}}\right\} 
\sinh \left\{ \sum\limits_{s=0}^{\infty}
\frac{\mathcal{E}_{2s+1}(a,z)
+d_{2s+1}(\alpha)}
{u^{2s+1}}\right\},
\end{equation}
as $u \to \infty$, uniformly for $0 \leq \arg(z) \leq \pi$ with $|z-z_{1}| \geq \delta >0$, under the condition \cref{eq12}. See \cite[Remark 1]{Dunster:2025:SAR} on how these expansions can be extended to $|z-z_{1}| \leq \delta$.

From \cite[Lemma 3.1]{Dunster:2025:SAR} we have an important result which is required in the next section:
\begin{lemma}
\label{lem:mero}
Each $(z-z_{1})^{1/2}\{\mathrm{E}_{2s+1}(\alpha,\phi)+d_{2s+1}(\alpha)\}$ ($s=0,1,2,\ldots$), regarded as a function of $z$, is meromorphic at $z=z_{1}$.
\end{lemma}

\section{Uniform asymptotic expansions for the zeros}
\label{sec:zeros}
We derive uniform asymptotic expansions for the complex zeros of $\theta_{n}(z;a)$ for large $n$, which are uniformly valid for unrestricted $z$ subject to \cref{eq15}.
To do so, we use the method of \cite{Dunster:2024:AZB}. We begin by defining a function $\mathcal{Z}(u,a,z)$ and coefficient $l_{n}(a)$ by the pair of equations
\begin{equation}
\label{eq37}
(-1)^{n} \frac{e^{a\pi i}}{n!} w^{(0)}_{n}(uz; a)
=e^{-\pi i/3} l_{n}(a)
\left\{\pdv{\mathcal{Z}(u,a,z)}{z}\right\}^{-1/2}
\mathrm{Ai}\left(u^{2/3}
\mathcal{Z}(u,a,z)\right),
\end{equation}
\begin{equation} 
\label{eq38}
\frac{1}{\Gamma(n + a - 1)} w^{(1)}_{n}(uz; a)
=l_{n}(a)\left\{\pdv{\mathcal{Z}(u,a,z)}{z}\right\}^{-1/2}
\mathrm{Ai}_{1}\left(u^{2/3}
\mathcal{Z}(u,a,z)\right).
\end{equation}
The factor $\{\partial \mathcal{Z}(u,a,z)/\partial z \}^{-1/2}$, along with Airy's equation \cite[Eq.~9.2.1]{NIST:DLMF}, ensures that both functions on the RHS of this pair of equations satisfy a linear second-order differential equation (with independent variable $z$) that has no first derivative term, which matches the same property of the functions on the LHS of these equations, namely the differential equation \cref{eq10}.

Now, from (\cite[Eq.~9.2.12]{NIST:DLMF}),
\begin{equation}
\label{eq39}
\mathrm{Ai}(z) 
=e^{\pi i/3}\mathrm{Ai}_{1}(z) 
+e^{-\pi i/3}\mathrm{Ai}_{-1}(z),
\end{equation}
and hence from the connection formula \cite[Eq.~(2.13)]{Dunster:2001:GBP}
\begin{equation} 
\label{eq40}
w^{(-1)}_{n}(z; a) = (-1)^{n+1} \frac{e^{a\pi i}}{n!} 
w^{(0)}_{n}(z; a) 
+ \frac{1}{\Gamma(n + a - 1)} w^{(1)}_{n}(z; a),
\end{equation}
we have
\begin{equation} 
\label{eq41}
w^{(-1)}_{n}(uz; a) = e^{\pi i/3} l_{n}(a)
\left\{\pdv{\mathcal{Z}(u,a,z)}{z}\right\}^{-1/2}
\mathrm{Ai}_{-1}\left(u^{2/3}
\mathcal{Z}(u,a,z)\right).
\end{equation}

Moreover, $\zeta \to \infty$ as $z \to 0$ or $z \to \pm\infty$; in these cases $Z \sim \zeta$ (see \cite{Dunster:2024:AZB}). The fundamental property of \cref{eq37,eq38,eq40} is that both functions in each identity are recessive at the same singularities, namely $z=-\infty,+\infty, 0$, respectively.

In studying the zeros we do not require the constant $l_{n}(a)$, but for completeness we note that it can be determined as follows. From \cref{eq03}, \cref{eq06}, and \cite[Eq.~13.2.34]{NIST:DLMF}
\begin{equation} 
\label{eq42}
\mathscr{W}\left\{w_{n}^{(0)}(z;a),w_{n}^{(-1)}(z;a)\right\}
=\frac{2^{1-2n-a}}{\Gamma(n+a-1)}.
\end{equation}
Thus from \cref{eq37}, \cref{eq41}, \cref{eq42} and \cite[Eq.~9.2.8]{NIST:DLMF}
\begin{equation} 
\label{eq43}
l_{n}(a)=\frac{e^{5\pi i/6} e^{(u+a)\pi i/2}u^{1/6}}
{2^{n-1}}\left\{\frac{\pi}
{2^{a} n!\Gamma(n+a-1)}\right\}^{1/2}.
\end{equation}

We focus on \cref{eq37}. Now, from \cite[Thm. 2.2]{Dunster:2024:AZB}, we have
\begin{equation}
\label{eq44}
\mathcal{Z}(u,a,z) \sim \zeta
+\sum_{s=1}^{\infty}
\frac{\Upsilon_{s}(a,z)}{u^{2s}}
\quad (u \to \infty),
\end{equation}
where each $\Upsilon_{s}(a,z)$ is analytic at the turning point $z=z_{1}$. This expansion is uniformly valid in an unbounded domain that includes the upper half-plane $0 \leq \arg(z) \leq \pi$. Thus, we shall use it to approximate all the zeros of $\theta_{n}(z;a)$ with nonnegative imaginary part, with those lying in the lower half-plane simply being the conjugates of these.

The coefficients $\Upsilon_{s}(a,z)$ are given by \cite[Thm. 2.2]{Dunster:2024:AZB}, with $\hat{E}_{2s+1}(z)$ replaced by $\mathrm{E}_{2s+1}(\alpha,\phi)+d_{2s+1}(\alpha)$ ($s=0,1,2,\ldots$). The inclusion of the constants $d_{2s+1}(\alpha)$ is required to ensure the required property stated in \cref{lem:mero}. From \cite[Eqs.~(2.28), (2.40) - (2.42)]{Dunster:2024:AZB} the first four coefficients in the series in $\cref{eq44}$ are given by
\begin{equation}
\label{eq45}
\Upsilon_{1}=
\frac{3 \xi (\mathrm{E}_{1}+d_{1})}
{2\zeta^{2}}
-\frac{5}{48\zeta^{2}},
\end{equation}
\begin{equation}
\label{eq46}
\Upsilon_{2}=
-\frac{\Upsilon_{1}^{2}}{4 \zeta}
+\frac{5 \Upsilon_{1}}{32\zeta^{3}}
+\frac{3 \xi (\mathrm{E}_{3}+d_{3})}{2\zeta^{2}}
-\frac{1105}{9216\zeta^{5}},
\end{equation}
\begin{equation}
\label{eq47}
\Upsilon_{3}=
-\frac{\Upsilon_{1}\Upsilon_{2}}{2 \zeta}
+\frac{\Upsilon_{1}^{3}}{24 \zeta^{2}}
-\frac{25\Upsilon_{1}^{2}}{128 \zeta^{4}}
+\frac{5 \Upsilon_{2}}{32\zeta^{3}}
+\frac{1105 \Upsilon_{1}}{2048 \zeta^{6}}
+\frac{3 \xi(\mathrm{E}_{5}+d_{5})}{2 \zeta^2}
-\frac{82825}{98304 \zeta^{8}},
\end{equation}
and
\begin{multline}
\label{eq48}
\Upsilon_{4}=
-\frac{\Upsilon_{1}^{4}}{64 \zeta^{3}}
+\frac{\Upsilon_{1}^{2}
\Upsilon_{2}}{8\zeta^{2}}
+\frac{175 \Upsilon_{1}^{3}}{768 \zeta^{5}}
-\frac{\Upsilon_{1}\Upsilon_{3}}{2\zeta}
-\frac{25\Upsilon_{1}\Upsilon_{2}}{64\zeta^{4}}
-\frac{\Upsilon_{2}^{2}}{4\zeta}
\\
-\frac{12155 \Upsilon_{1}^{2}}
{8192 \zeta^{7}}
+\frac{5\Upsilon_{3}}{32\zeta^{3}}
+\frac{1105\Upsilon_{2}}{2048\zeta^{6}}
+\frac{414125\Upsilon_{1}}{65536\zeta^{9}}
+\frac{3\xi(\mathrm{E}_{7}+d_{7})}{2\zeta^{2}}
-\frac{1282031525}{88080384\zeta^{11}},
\end{multline}
with subsequent ones given by \cite[Thm. 2.2]{Dunster:2024:AZB}.

Next, let $t_{m}(u,a)$ ($m=1,2,3, \ldots ,\lfloor (n+1)/2\rfloor$) be the complex zeros of $\theta_{n}(t;a)$ in the upper half-plane $\Re(t) \geq 0$, so that
\begin{equation}
\label{eq49}
\theta_{n}(t_{m}(u,a);a)=0
\quad (m=1,2,3,\dots ,\lfloor (n+1)/2\rfloor).
\end{equation}
From \cref{eq03,eq37} they satisfy the implicit equation
\begin{equation}
\label{eq50}
\mathcal{Z}(u,a,u^{-1}t_{m}(u,a))
=u^{-2/3}\mathrm{a}_{m}
\quad (m=1,2,3,\ldots),
\end{equation}
where $x=\mathrm{a}_{m}$ is the $m$th negative zero of $\mathrm{Ai}(x)$ ordered by increasing absolute values.

From \cite[Thm. 3.1]{Dunster:2024:AZB} we obtain the uniform asymptotic expansion we seek, namely
\begin{equation}
\label{eq51}
t_{m}(u,a) \sim
u \sum_{s=0}^{\infty}
\frac{\tau_{m,s}(\alpha)}{u^{2s}}
\quad (u \to \infty, \, m=1,2,3,\ldots),
\end{equation}
for coefficients $\tau_{m,s}(\alpha)$ which we determine next.

With the branch cut for $Z$ as described in \cref{sec:LG}, and principal branches for the logarithms we have $\xi=-2i|\mathrm{a}_{m}|^{3/2}/(3u)$ for $\zeta=u^{-2/3}\mathrm{a}_{m}$. Thus, on plugging \cref{eq51} into \cref{eq50}, using \cref{eq16,eq17,eq44}, re-expanding in inverse powers of $u$, and then equating like powers, we can find in turn the coefficients for each prescribed $m \in \{1,2,3,\ldots \lfloor(n+1)/2\rfloor \}$.

Consequently, the leading term $\tau_{m,0}=\tau_{m,0}(\alpha)$ is given implicitly by
\begin{multline}
\label{eq52}
-\frac{2 i}{3 u} |\mathrm{a}_{m}|^{3/2}
=Z_{m,0}+\left(1+\frac12 \alpha\right)
\ln\left\{\frac{\tau_{m,0}}
{4Z_{m,0}+2\alpha(Z_{m,0}+\tau_{m,0}+2)+4
+\alpha^{2}} \right\} 
\\
+\frac12 \alpha\left\{ \ln\left(-2Z_{m,0}
-2\tau_{m,0}-\alpha\right)+\pi i\right\}
+\frac{1}{2}\ln(1+\alpha)
+\left(2+\frac12 \alpha\right)\ln(2)
-\frac{1}{2}(1+\alpha)\pi i,
\end{multline}
in which $Z_{m,0}=Z(\tau_{m,0})$, where $Z(z)$ is given by \cref{eq17} with the branch as described below that equation. Thus
\begin{equation} 
\label{eq53}
Z_{m,0}=Z_{m,0}(\alpha) = -\left\{(\tau_{m,0} - z_1)
(\tau_{m,0} - z_2)\right\}^{1/2},
\end{equation}
which is negative for $-\infty<\tau_{m,0}<0$, since $z=\tau_{m,0}$ lies to the left of the cut along the anti-Stokes line $\Im(\xi)=0$ from $z=z_{1}$ to $z=0$. For numerical purposes, in \cref{eq52} the second logarithm was expressed in a form that ensures the correct branch for these complex values of $\tau_{m,0}$.

Next, let
\begin{equation}
\label{eq54}
\zeta_{m,0}=\zeta\left(\tau_{m,0}\right)
=u^{-2/3}\mathrm{a}_{m}, \;
\zeta_{m,0}^{\prime}
=\zeta'\left(\tau_{m,0}\right), \;
\zeta_{m,0}^{\prime \prime}
=\zeta''\left(\tau_{m,0}\right), \;
\ldots,
\end{equation}
and similarly for $s=1,2,3,\ldots$ let
\begin{equation}
\label{eq55}
\Upsilon_{m,s}
=\Upsilon_{s}\left(\tau_{m,0}\right), \;
\Upsilon_{m,s}^{\prime}
=\Upsilon_{s}'\left(\tau_{m,0}\right), \;
\Upsilon_{m,s}^{\prime \prime}
=\Upsilon_{s}''\left(\tau_{m,0}\right), \;
\ldots.
\end{equation}
Then we can apply \cite[Thm. 3.1]{Dunster:2024:AZB} to find the coefficients in \cref{eq51}. Having determined $\tau_{m,0}$, the next four coefficients are of the same form as \cite[Eqs.(3.47) - (3.50)]{Dunster:2024:AZB} (with the appropriate change of notation). Thus
\begin{equation}
\label{eq56}
\tau_{m,1}=
-\frac{\Upsilon_{m,1}}
{\zeta'_{m,0}},
\end{equation}
\begin{equation}
\label{eq57}
\tau_{m,2} =-\frac{1}
{2 \zeta'_{m,0}}
\left\{
\tau_{m,1}^{2} \, 
\zeta''_{m,0}+2 \tau_{m,1}\, \Upsilon'_{m,1}
+2\Upsilon_{m,2}
\right\},
\end{equation}
\begin{multline}
\label{eq58}
\tau_{m,3} =-\frac{1}
{6 \zeta'_{m,0}}
\left\{
\tau_{m,1}^3 \,\zeta'''_{m,0} + 6\tau_{m,1}\,\tau_{m,2}\,
\zeta''_{m,0} + 3\tau_{m,1}^2\Upsilon''_{m,1}
\right.
\\
\left. 
+ 6\tau_{m,2}\Upsilon'_{m,1} 
+ 6\tau_{m,1}\Upsilon'_{m,2} 
+ 6\Upsilon_{m,3} 
\right\},
\end{multline}
and\footnote{There is a misprint in \cite[Eq.~(3.50)]{Dunster:2024:AZB}}
\begin{multline}
\label{eq59}
\tau_{m,4} =-\frac{1}
{24 \zeta'_{m,0}}
\left\{
\tau_{m,1}^4 \,\zeta^{(4)}_{m,0}
+12 \tau_{m,1}^{2}\, \tau_{m,2}\,\zeta'''_{m,0}
+24\tau_{m,1}\, \tau_{m,3}\,\zeta''_{m,0}
\right.
\\
+12\tau_{m,2}^2\,\zeta''_{m,0}
+4\tau_{m,1}^3 \,\Upsilon'''_{m,1}
+24\tau_{m,1}\,\tau_{m,2}\,\Upsilon''_{m,1}
+12\tau_{m,1}^2 \,\Upsilon''_{m,2}
\\
\bigl. 
+24\tau_{m,3}\,\Upsilon'_{m,1}
+24\tau_{m,2}\,\Upsilon'_{m,2}
+24\tau_{m,1}\,\Upsilon'_{m,3}
+24\Upsilon_{m,4}
\Bigr\}.
\end{multline}
These, of course, result in different coefficients than in the Bessel function case, due to the difference here in the variable $\zeta$ as well as the $\Upsilon$ coefficients.

\subsection{Numerical examples}
\label{sec:numerics}
We now approximate $t_{m}(u,a)$ by the series \cref{eq51}. In order to do so, we require the derivatives with respect to $z$ of $\zeta$, $\xi$, $\phi$ and $\Upsilon_{s}$. For the latter, it is convenient to let $\rho=1/\zeta$. Then from \cref{eq16}
\begin{equation} 
\label{eq60}
\rho'=-\tfrac32 \rho^{4}\xi\xi'.
\end{equation}
For the other derivatives, we use 
\begin{equation} 
\label{eq61}
\zeta'=\frac{2\xi' \zeta}{3\xi},
\end{equation}
\begin{equation} 
\label{eq62}
\xi'=f^{1/2}(a,z)=\frac{\sigma}{z\sin(\phi)},
\end{equation}
and 
\begin{equation} 
\label{eq63}
\phi'=-\frac{\sin^{2}(\phi)}{\sigma},
\end{equation}
which follow from \cref{eq11,eq14,eq16,eq18,eq20}. 

We compute the first five terms $\tau_{m,s}$ ($s=0,1,2,3,4$) in the expansion \cref{eq51}. This is achieved by the following steps.

\begin{itemize}
 \item For each prescribed $n$ and $m \in \{1,2,3,\ldots,\lfloor n+1\rfloor\}$ find $\xi=-2i|\mathrm{a}_{m}|^{3/2}/(3u)$ where $u=n+\frac12$, and $\zeta=\zeta_{m,0}=\mathrm{a}_{m}u^{-2/3}$. Use these values for $\xi$ and $\zeta$ below.
 \item Use \cref{eq52} to numerically evaluate $\tau_{m,0}$. To obtain the correct root we found it efficient to set $\tau_{m,0}=-0.5+w$ in the equation and then numerically solve for $w$. Then use $z=\tau_{m,0}$ in what follows.
 \item Use \cref{eq17,eq18,eq19,eq14,eq61,eq62,eq63} to compute $\zeta'_{m,0}$, $\zeta''_{m,0}$, $\zeta'''_{m,0}$ and $\zeta^{(4)}_{m,0}$.
 \item From \cref{eq21,eq22,eq23,eq24,eq63} for $z=\tau_{m,0}$ compute $\mathrm{E}_{1}$ and its first three $z$ derivatives, $\mathrm{E}_{3}$ and its first two $z$ derivatives, $\mathrm{E}_{5}$ and its $z$ derivative, and $\mathrm{E}_{7}$.
 \item Use \cref{eq27,eq28,eq30,eq31,eq32,eq33,eq45,eq46,eq47,eq48,eq60,eq62,eq63} and the above values to compute in turn $\Upsilon_{m,1}$, $\Upsilon'_{m,1}$, $\Upsilon''_{m,1}$, $\Upsilon'''_{m,1}$, $\Upsilon_{m,2}$, $\Upsilon'_{m,2}$, $\Upsilon''_{m,2}$, $\Upsilon_{m,3}$, $\Upsilon'_{m,3}$, and $\Upsilon_{m,4}$. 
 \item Evaluate ${\tau}_{m,s}$ for $s=1,2,3,4$ in turn from \cref{eq56,eq57,eq58,eq59}. 
\end{itemize}

\vspace*{0.2cm}

Numerical examples of the approximations $t_m$ obtained using the scheme described above are presented in 
\cref{table:1,table:2} for $a=1.01$, $20.2$, respectively, and different values of $n$ and $m$. The approximations, 
implemented in Maple\footnote{The file can be obtained from {https://github.com/AmparoGil/AsympZerosRGBPs}}, are compared with the values obtained using the numerical algorithm described in \cite{Segura:2013:CCZ} (also implemented in Maple). The relative errors from these comparisons are shown in the tables. Relative errors are all close to or better than $10^{-15}$.

\begin{table}
$$
\begin{array}{lcl}
\hline
n,\,m & t_m & \mbox{Relative error} \\
 \hline
15,\,1 & -3.1559515225814951808+12.586271690843017387i & 1.8 \times 10^{-15} \\
15,\,3 & -6.9360218173803455640+8.6292759166638006520i & 6.1\times 10^{-16}\\
\hline
30,\,1 & -4.2425750716206130472+27.006358468998877565i & 1.1\times 10^{-16} \\
30,\,3 & -9.7584463264409865096+22.392832031435945931i & 1.3 \times 10^{-16} \\
30,\,10 & -18.102790325129739597+9.4722422021510892034i & 7.7 \times 10^{-17} \\
30,\,15 & -19.702854218331257062+.85611271550820061202i & 2.9 \times 10^{-16} \\
\hline
50,\,1 &-5.2055266715795128190+46.482961682470093754i & 1.9 \times 10^{-17} \\
50,\,3 & -12.181102558122645217+41.239145916888100131i & 7.9 \times 10^{-17} \\
50,\,10 & -24.683402130958153499+27.225504025486397962i & 1.4 \times 10^{-16} \\
50,\,15 & -29.379559025204265717+18.222895815367965462i & 2.2 \times 10^{-16} \\
50,\,25 & -32.962750529211803345+.86074820845854851940i & 6.7 \times 10^{-18} \\
\end{array}
$$
\caption{Approximations to the zeros of the reverse generalized Bessel polynomials for $a=1.01$ and different values of $n$ and $m$.}
\label{table:1}
\end{table}

\begin{table}
$$
\begin{array}{lcl}
\hline
n,\,m & t_m & \mbox{Relative error} \\
 \hline
15,\,1 & -12.715856054909203812+18.788546633810651464i & 1.1 \times 10^{-15} \\
15,\,3 &-16.514653825298059143+12.612556755577648289i & 2.6\times 10^{-15}\\
\hline
30,\,1 & -13.800334806578766149+34.380365451162645216i & 3.4\times 10^{-16} \\
30,\,3 & -19.310221900147056579+28.210989284732813206i & 2.8 \times 10^{-16} \\
30,\,10 & -27.717880396627235555+11.750965665786499280i & 5.4 \times 10^{-17} \\
30,\,15 & -29.339399892921113584+1.0590134228243351098i & 2.8 \times 10^{-16} \\
\hline
50,\,1 & -14.766307319696546646+54.504885286408130512i
 & 2.9 \times 10^{-16} \\
50,\,3 & -21.724567399352576652+48.087744580616150218i
 & 6.1 \times 10^{-17} \\
50,\,10 & -34.260698846474016613+31.438165321383787957i
 & 1.6 \times 10^{-16} \\
50,\,15 & -38.989834370513922989+20.967450446744804559i
 & 2.3 \times 10^{-16} \\
50,\,25 & -42.605131456252572254+.98772884689217274567i
 & 2.6 \times 10^{-18} \\
\end{array}
$$
\caption{Approximations to the zeros of the reverse generalized Bessel polynomials for $a=20.2$ and different values of $n$ and $m$.}
\label{table:2}
\end{table}

In the implementation of the numerical algorithm, the reverse generalized Bessel polynomials are evaluated through their relation with the generalized Laguerre polynomials \cite[\S~18.34(i)]{NIST:DLMF}
\begin{equation}
\label{eq64}
\theta_{n}(z;a)=\left(-\tfrac{1}{2}\right)^{n}
n!L_n^{(1-2n-a)}(2z).
\end{equation}

\begin{figure}
 \centering
 \includegraphics[width=0.8\textwidth,keepaspectratio]{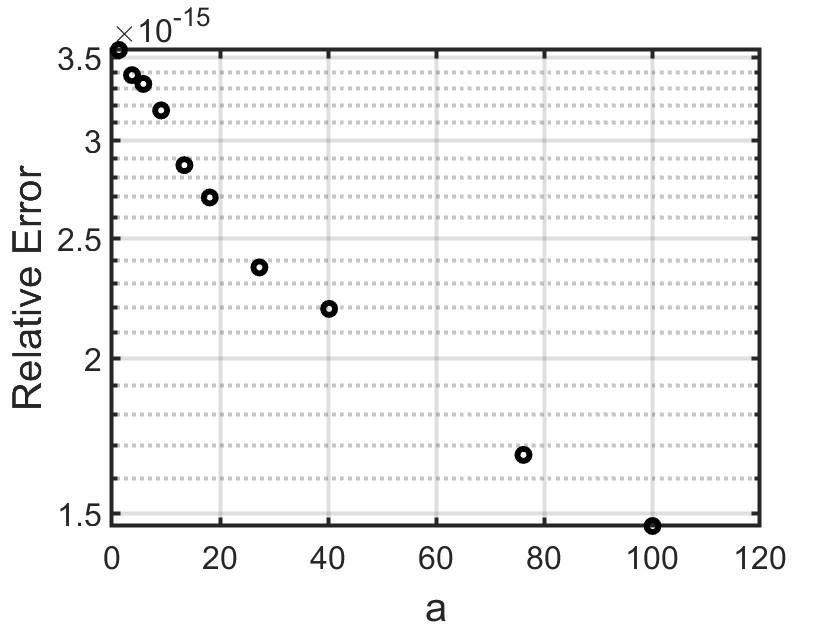}
 \caption{ Relative errors as a function of $a$ (with $n$ fixed at $15$ and $m$ at $3$).}
 \label{fig:fig1}
\end{figure}

Furthermore, \cref{fig:fig1} displays the relative error of the approximations to the zeros for various values of $a$ (with $n$ fixed at $15$ and $m$ at $3$). It is observed that the relative error decreases as $a$ increases, and the maximum relative error attained is less than $4 \times 10^{-15}$.

\begin{figure}
 \centering
 \includegraphics[width=0.99\textwidth,keepaspectratio]{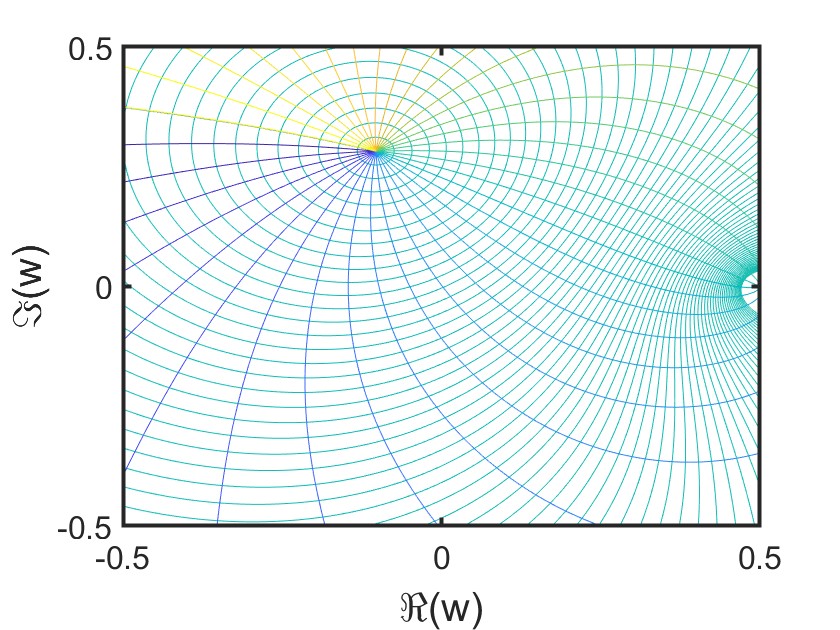}
 \caption{ Plot of the function $F(w)$ for $a=1.01$, $m=10$ and $n=30$. }
 \label{fig:fig2}
\end{figure}

Notably, equation \cref{eq52} can be solved using Newton's method to obtain the corresponding root. By expressing the equation as $F(w)=0$, \cref{fig:fig2} illustrates, as an example, a plot of the function $F(w)$ for $a=1.01$, $m=10$ and $n=30$. The root obtained by applying Newton's method is $w=-0.0935299175+0.310545771i$, which agrees with the location of the zero of $F(w)$ shown in the figure.

\section{A numerical algorithm to compute the zeros of \texorpdfstring{$\theta_n(z;a)$}{theta\_n(z;a)}}
\label{sec:algorithm}

Given that the asymptotic approximations obtained exhibit high accuracy for moderate to large values of $n$, they can be employed as initial estimates in the numerical algorithm described in \cite{Segura:2013:CCZ} devised to compute all the zeros. 
Within this algorithm, the fundamental components are an iteration function $T_{n}(a,z)$ and a step function $H_{n}(a,z)$:
\begin{equation}
\label{eq65}
T_{n}(a,z)=z-\frac{1}{\sqrt{\Omega_{n}(a,z)}}
\arctan\left(\frac{\sqrt{\Omega_{n}(a,z)}\,w_{n}^{(0)}(z;a)}
{\partial w_{n}^{(0)}(z;a)/\partial z}\right),
\end{equation}
\begin{equation}
\label{eq66}
H_{n}(a,z)=z+\frac{\pi}{\sqrt{\Omega_{n}(a,z)}},
\end{equation}
where
\begin{equation}
\label{eq67}
\Omega_{n}(a,z)=-1+\frac{2-a}{z}-\left(n+\frac{1}{2}a \right)
\left(n+ \frac{1}{2}a-1 \right)\frac{1}{z^2}.
\end{equation}

The function $w_{n}^{(0)}(z;a)$ appearing in \cref{eq65} is related to the reverse generalized Bessel polynomials $\theta_n(z;a)$ by \cref{eq03}. It is therefore clear that $w_{n}^{(0)}(z;a)$ can be used to locate the zeros of $\theta_n(z;a)$. In \cite{Dunster:2021:CPB}, we discussed the possible use of asymptotic expansions to compute the reverse generalized Bessel polynomials appearing in the iteration function \cref{eq65}. However, this is not necessary if accurate approximations are available for the first zero (combined with the use of Taylor series). This is precisely the strategy we discuss next. The resulting algorithm, which avoids the explicit computation of the function in order to locate its zeros, is both accurate and highly efficient.

We use as a starting point for finding the Taylor series the equation satisfied by $w_{n}=w_{n}^{(0)}(z;a)$, written in the form
\begin{equation}
Pw_{n}''+Q_{n}w_{n}=0,
\label{eq69}
\end{equation}
where 
\begin{equation}
P=P(z)=z^2,\;
Q_{n}=Q_{n}(a,z)=z^2\Omega_{n}(a,z).
\label{eq70}
\end{equation}

We have $P^{(m)}=Q_{n}^{(m)}=0$ for $m>2$. Then, we obtain for $w_{n}=w_{n}^{(0)}(z;a)$
\begin{equation}
\label{eq71}
\sum_{m=0}^{2}\binom{j}{m} 
\left\{P^{(m)}w_{n}^{(j+2-m )} 
+ Q_{n}^{(m)}w_{n}^{(j-m)}\right\}=0,
\end{equation}
which gives the following recurrence relation for the derivatives
\begin{multline}
\label{eq72}
z^2w_{n}^{(k+2)}+2kzw_{n}^{(k+1)}
+\left\{Q_{n}+k(k-1)\right\}w_{n}^{(k)}
\\
-k(2z+a-2) w_{n}^{(k-1)} -k(k-1)w_{n}^{(k-2)}=0
\quad (k=2,3,4,\ldots).
\end{multline}

This recurrence relation, which is not ill-conditioned as $k$
becomes large, can be used to compute the derivatives appearing in the Taylor series
\begin{multline}
\label{eq73}
w_{n}(z_{i+1})=\sum_{k=0}^{N}w_{n}^{(k)}(z_i)\Frac{h^k}{k!}
+ {\cal O} (h^{N+1}),\\
w_{n}'(z_{i+1})=\sum_{k=0}^{N}w_{n}^{(k+1)}(z_i)\Frac{h^k}{k!}
+ {\cal O} (h^{N+1}).
\end{multline}

To apply \cref{eq72} in the algorithm step by step, we need
$w_{n}(z_i)$ and $w_{n}'(z_i)$, which are known from the previous step, and (again suppressing the $a$-dependence)
\begin{multline}
\label{eq74}
w_{n}''(z_i)=-\Frac{Q_{n}(z_i)}{P(z_i)}w_{n}(z_i),\\
w_{n}'''(z_i)=-\Frac{Q_{n}(z_i)}{P(z_i)}w_{n}'(z_i)
+\left\{\Frac{(2-a)}{z_i^2} -\Frac{2\left(n+\frac12 a\right)\left(n+\frac12 a-1\right)}
{z_i^3}\right\}w_{n}(z_i).
\end{multline}

The resulting algorithm works as follows: After a zero $z_i$ has been obtained (using an asymptotic approximation for the first zero), the next zero $z_{i+1}$ is obtained by taking the step $H_{n}(a,z)$ given in \cref{eq66} and then iterating $T_{n}(a,z)$ given by \cref{eq65} until convergence is reached. Taylor series \cref{eq73} are used to calculate the functions $w_{n}(z)$ and $w_{n}'(z)$ needed in the iteration function, and the values $w_{n}(z_i)=0$ and $w_{n}'(z_i)=1$ are used in this computation. A Matlab implementation of the algorithm can be obtained from GitHub\footnote{{https://github.com/AmparoGil/NumerZerosRGBPs}}. In the present implementation, and for reasons of computational efficiency, only the first three coefficients of the uniform expansion \cref{eq51} are employed in order to approximate the first zero to be computed (namely, the one with the largest real part in absolute value). The algorithm computes the complex zeros with nonnegative imaginary part (as it was mentioned before, the remaining complex zeros are simply the complex conjugates of those obtained).

\begin{figure}
 \centering
 \includegraphics[width=0.8\textwidth,keepaspectratio]{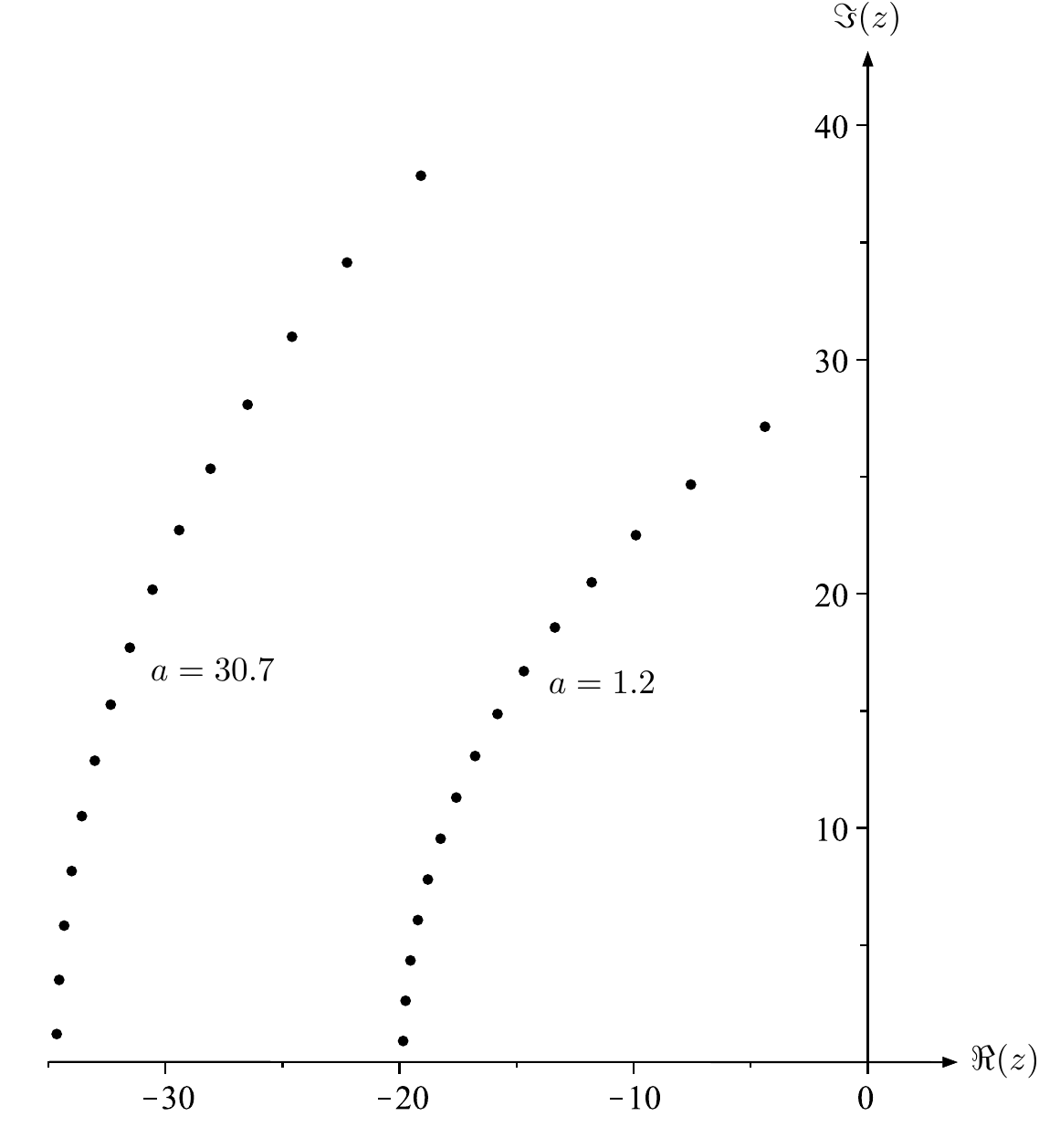}
 \caption{Zeros in the second quadrant obtained by the numerical algorithm for $n=30$ and $a=1.2, \, 30.7$.}
 \label{fig:fig3a}
\end{figure}

\begin{figure}
 \centering
 \includegraphics[width=0.8\textwidth,keepaspectratio]{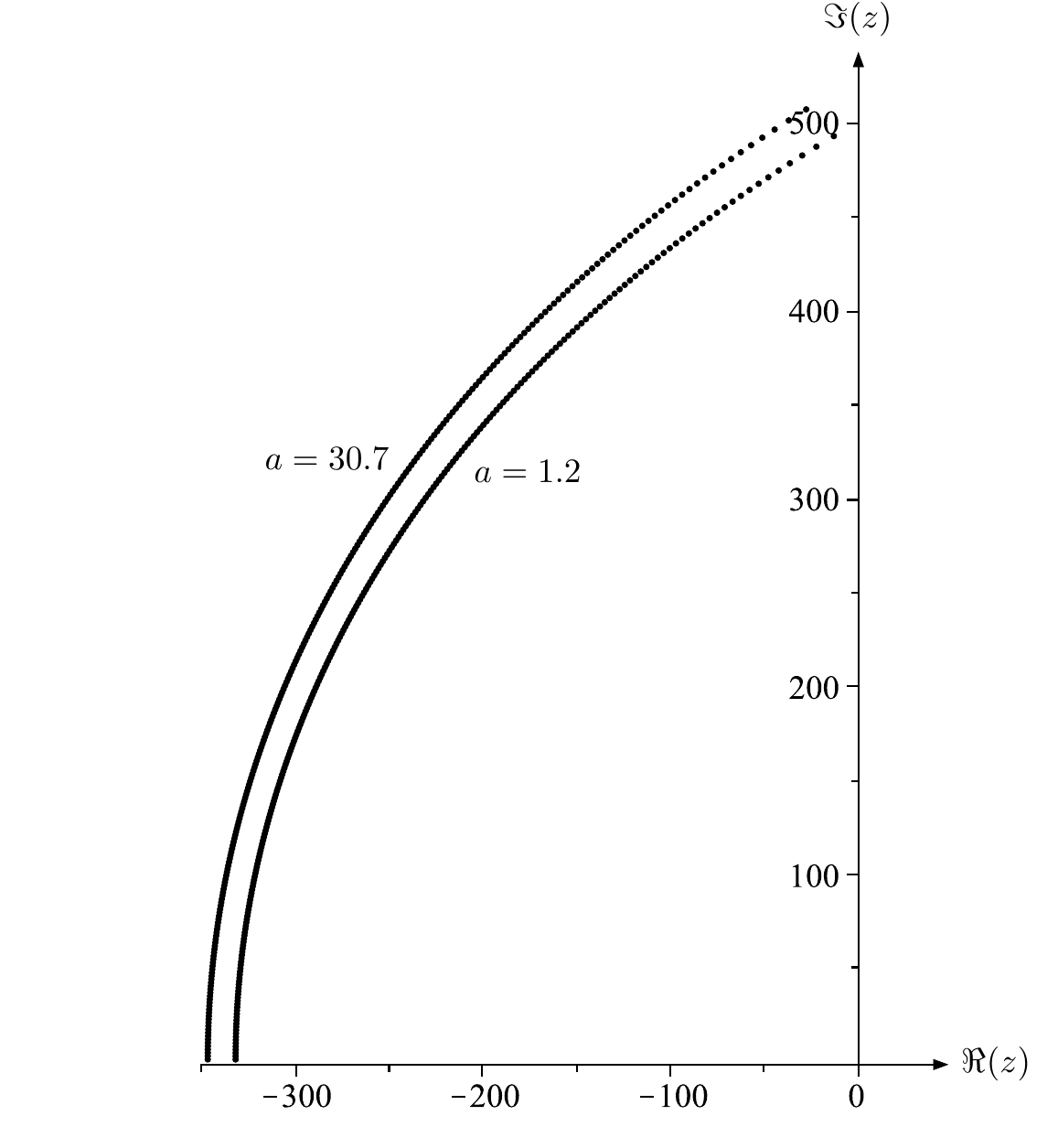}
 \caption{Zeros in the second quadrant obtained by the numerical algorithm for $n=500$ and $a=1.2, \, 30.7$.}
 \label{fig:fig3b}
\end{figure}

\begin{figure}
 \centering
 \includegraphics[width=\textwidth,keepaspectratio]{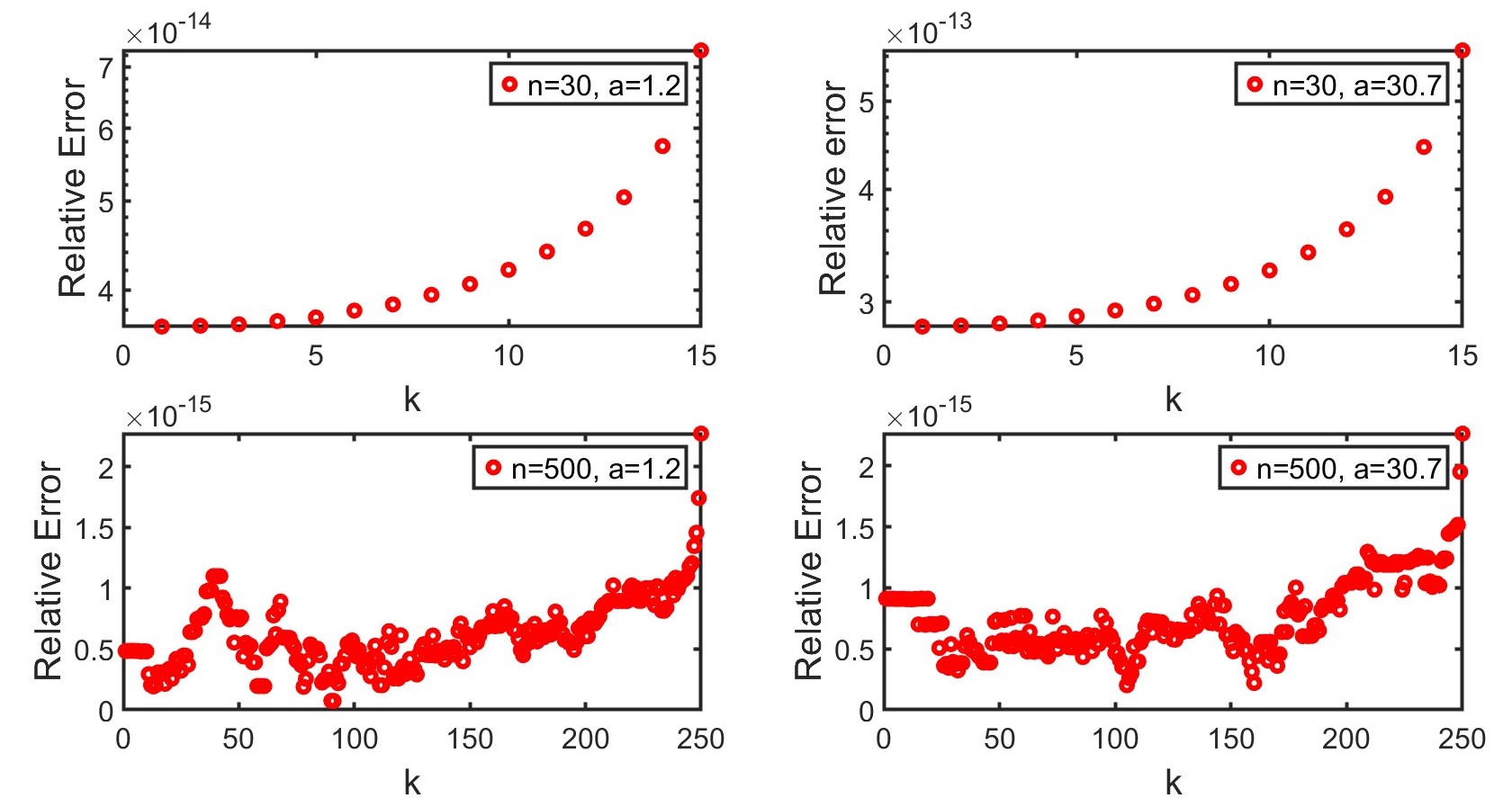}
 \caption{ Relative errors of the computed zeros obtained using the numerical algorithm. }
 \label{fig:fig4}
\end{figure}

Plots of the zeros of $\theta_{n}(z;a)$ in the second quadrant obtained by the numerical algorithm are shown in \cref{fig:fig3a,fig:fig3b} for $n=30$, $n=500$, $a=1.2$ and $a=30.7$. The accuracy prescribed for the iterations with $T_{n}(a,z)$ is $\epsilon=10^{-12}$, although the precision attained for the zeros is higher when the value of $n$ is sufficiently large, as can be observed in \cref{fig:fig4}. The relative errors have been obtained by comparison with the Maple implementation of the numerical algorithm used in \cref{sec:numerics}. The number of digits was set to $60$ in the Maple implementation in order to ensure the correct
evaluation of the quotients appearing in the iteration function. It can be observed that, even when a large number of zeros are computed, the relative errors remain well controlled.

Finally, \cref{table:3} presents representative CPU times obtained from the execution of our numerical algorithm for several values of $n$. CPU times were measured by executing our algorithm in Matlab R2024b on a Dell Latitude $7410$ equipped with $16$ GB RAM and an Intel Core i$5$-$10210$U CPU at $1.6$ GHz. Notably, increasing from computing $15$ zeros ($n=30$) to $1000$ zeros ($n=2000$) raises CPU time by only a factor of $10$, despite the number of zeros increasing by a far larger factor. This observation provides clear evidence of the efficiency of the algorithm developed for the computation of the zeros.

\begin{table}
$$
\begin{array}{cl}
\hline
n & \mbox{CPU time ($s$)} \\
 \hline
30 & 2.8 \times 10^{-3} \\
200 & 5.9 \times 10^{-3} \\
500 & 1.1 \times 10^{-2} \\
1000 & 1.9 \times 10^{-2} \\
2000 & 3.7 \times 10^{-2} \\
\end{array}
$$
\caption{Representative CPU times obtained from the execution of our numerical algorithm for $a=2.3$ and several values of $n$.}
\label{table:3}
\end{table}

\FloatBarrier

\section*{Acknowledgment}
Financial support from Ministerio de Ciencia e Innovación pro\-jects PID2021-127252NB-I00 
and PID2024-159583NB-I00\\
(MICIU/ AEI / 10.13039/501100011033 / FEDER, UE) are acknowledged.

\section*{Conflict of interest}
The authors have no conflict of interest to declare that is relevant to the content of this article.

\FloatBarrier

\makeatletter
\interlinepenalty=10000
\bibliographystyle{siamplain}
\bibliography{biblio}

\begin{thebibliography}{10}

\bibitem{Carimalo:2018:MFG}
{\sc C.~Carimalo}, {\em Maximally flat group delay of {B}essel polynomials}, Circ. Syst. Signal Pr., 37 (2018), pp.~5174--5177, \url{https://doi.org/10.1007/s00034-018-0812-x}.

\bibitem{Carpenter:1992:GBP}
{\sc A.~Carpenter}, {\em Asymptotics for the zeros of the generalized {B}essel polynomials}, Numer. Math., 62 (1992), p.~465–482, \url{https://doi.org/10.1007/BF01396239}.

\bibitem{deBruin:1981:ZGBI}
{\sc M.~G. de~Bruin, E.~B. Saff, and R.~S. Varga}, {\em On the zeros of generalized {B}essel polynomials. {I}}, Indag. Math., 84 (1981), pp.~1--13, \url{https://doi.org/10.1016/1385-7258(81)90013-5}.

\bibitem{deBruin:1981:ZGBII}
{\sc M.~G. de~Bruin, E.~B. Saff, and R.~S. Varga}, {\em On the zeros of generalized {B}essel polynomials. {II}}, Indag. Math., 84 (1981), pp.~14--33, \url{https://doi.org/10.1016/1385-7258(81)90014-7}.

\bibitem{NIST:DLMF}
{\em {\it NIST Digital Library of Mathematical Functions}}.
\newblock Release 1.1.6 of 2022-06-30, \url{http://dlmf.nist.gov/}.
\newblock F.~W.~J. Olver, A.~B. {Olde Daalhuis}, D.~W. Lozier, B.~I. Schneider, R.~F. Boisvert, C.~W. Clark, B.~R. Miller, B.~V. Saunders, H.~S. Cohl, and M.~A. McClain, eds.

\bibitem{Dunster:2001:GBP}
{\sc T.~M. Dunster}, {\em Uniform asymptotic expansions for the reverse generalized {B}essel polynomials, and related functions}, SIAM J. Math. Anal., 32 (2001), pp.~987--1013, \url{https://doi.org/10.1137/S0036141099359068}.

\bibitem{Dunster:2024:AZB}
{\sc T.~M. Dunster}, {\em Uniform asymptotic expansions for the zeros of {B}essel functions}, SIAM J. Math. Anal., 56 (2024), pp.~6521--6550, \url{https://doi.org/10.1137/23M1611269}.

\bibitem{Dunster:2025:SAR}
{\sc T.~M. Dunster}, {\em Simplified {A}iry function asymptotic expansions for reverse generalised {B}essel polynomials}, J. Classical Anal.,  (2025), \url{https://arxiv.org/abs/2506.20934}.
\newblock In Press.

\bibitem{Dunster:2021:CPB}
{\sc T.~M. Dunster, A.~Gil, D.~Ruiz-Antolin, and J.~Segura}, {\em Computation of the reverse generalized {B}essel polynomials and their zeros}, Comput. Math. Methods, 3(6):e1198 (2021), \url{https://doi.org/10.1002/cmm4.1198}.

\bibitem{Fila:2015:AFI}
{\sc I.~M. Filanovsky}, {\em Enhancing amplifiers/filters bandwidth by transfer function zeroes}, in 2015 IEEE International Symposium on Circuits and Systems (ISCAS), 2015, pp.~141--144, \url{https://doi.org/10.1109/ISCAS.2015.7168590}.

\bibitem{Johnson:1976:FIL}
{\sc J.~Johnson, D.~Johnson, P.~Boudra, and V.~Stokes}, {\em Filters using {B}essel-type polynomials}, IEEE Transactions on Circuits and Systems, 23 (1976), pp.~96--99, \url{https://doi.org/10.1109/TCS.1976.1084174}.

\bibitem{Kong:2024:EIG}
{\sc D.~Kong, J.~Shen, L.~Wang, and S.~Xiang}, {\em Eigenvalue analysis and applications of the {L}egendre dual-{P}etrov-{G}alerkin methods for initial value problems}, Adv. Comput. Math., 50 (2024), \url{https://doi.org/10.1007/s10444-024-10190-z}.

\bibitem{Martinez:1977:TFG}
{\sc J.~Martinez}, {\em Transfer functions of generalized {B}essel polynomials}, IEEE Trans. Circuits Syst., 24 (1977), pp.~325--328, \url{https://doi.org/10.1109/TCS.1977.1084347}.

\bibitem{Olver:1997:ASF}
{\sc F.~W.~J. Olver}, {\em Asymptotics and special functions}, AKP Classics, A K Peters Ltd., Wellesley, MA, 1997.
\newblock Reprint of the 1974 original [Academic Press, New York].

\bibitem{Pasquini:2000:GBP}
{\sc L.~Pasquini}, {\em Accurate computation of the zeros of the generalized {B}essel polynomials}, Numer. Math., 86 (2000), pp.~507--538, \url{https://doi.org/10.1007/s002110000166}.

\bibitem{Segura:2013:CCZ}
{\sc J.~Segura}, {\em Computing the complex zeros of special functions}, Numer. Math., 124 (2013), pp.~723--752, \url{https://doi.org/10.1007/s00211-013-0528-6}.

\end{thebibliography}
\end{document}